\documentclass[a4paper,11pt]{article}

\usepackage{enumerate}
\usepackage{amsmath}
\usepackage{amssymb}
\usepackage{amsfonts}
\usepackage{fancyhdr}

\newcommand{\cp}{\mathbb{C}}
\newcommand{\qp}{\mathbb{Q}}
\newcommand{\pp}{\mathbb{P}}
\newcommand{\np}{\mathbb{N}}

\newcommand{\zp}{\mathbb{Z}}

\newcommand{\wt}{\widetilde}

\newtheorem{lem}{Lemma}

\newtheorem{thr}{Theorem}
\newtheorem{rem}{Remark}

\newenvironment{preuve}
{\noindent{\textit{Proof. }}}
{\normalsize\rmfamily\\[.2cm]} 

\title{An interpolation theorem in toric varieties}
\author{WEIMANN Martin}

\begin{document}

\maketitle

\begin{abstract}
In the spirit of a theorem of Wood \cite{W:gnus}, we give necessary and sufficient conditions for a family of germs of analytic hypersurfaces in a smooth projective toric variety $X$ to be interpolated by an algebraic hypersurface with a fixed class in the Picard group of $X$. 
\end{abstract}

\section{Introduction}

Let $X$ be a compact algebraic variety over $\cp$. We are interested by the following problem:
\vskip 1mm
\noindent 
{\it Let $V_1,\ldots,V_N$ be a collection of germs of smooth analytic hypersurfaces at pairwise distincts smooth points $p_1,\ldots,p_N$ of $X$, and fix $\alpha$ in the Picard group $Pic(X)$ of $X$. When does there exist an algebraic hypersurface $\wt{V}\subset X$ with class $\alpha$ containing all the germs $V_i$ ?}
\vskip 1mm
\noindent
A natural way to approach this problem is to study sums and products of values of rational functions at points of intersection of the germs $V_i$ with a  "moving" algebraic curve\footnote{This philosophy has been initiated by Abel in his studies of abelian integrals \cite{ab:gnus}.}. 
\vskip 1mm
\noindent
Let us recall a theorem of Wood \cite{W:gnus} treating the case of $N$ germs in an affine chart $\cp^n$ of $X=\pp^n$, transversal to the line $l_0=\{x_1=\cdots=x_{n-1}=0\}$. Any line $l_a$ close to $l_0$ has affine equations $x_k = a_{k0}+a_{k1}x_{n}$, $k=1,\ldots,n-1$.
The trace on $V=V_1\cup\cdots\cup V_N$ of any function $f$ holomorphic in an analytic neighborhood of $V$ is the function
$$
a \longmapsto  Tr_V(f) (a):=\sum_{p\in V\cap l_a} f(p)\,,
$$
defined and holomorphic for $a=((a_{10},a_{11}),\ldots,(a_{n-1,0},a_{n-1,1}))$ close enough to $0\in \cp^{2n-2}$. 
\begin{thr} (Wood, \cite{W:gnus})
There exists an algebraic hypersurface $\wt{V}\subset \pp^n$ of degree $N$ which contains $V$ if and only if the function $a\mapsto  Tr_V(x_n)(a)$ is affine in the constant coefficients $a_{0}=(a_{10},\ldots,a_{n-1,0})$.
\end{thr}

We show here that Theorem 1 has a natural generalization to germs $V_1,\ldots,V_N$ in general position in a smooth toric compactification $X$ of $\cp^n$ endowed with an ample line bundle. As in \cite{W:gnus}, our proof gives an explicit construction of the polynomial equation of the interpolating hypersurface in the affine chart $\cp^n$. Moreover, we  characterize the class of $\wt{V}$ in $Pic(X)$.
\vskip 1mm
\noindent
On any projective variety $X$, there exist very ample line bundles $L_1,\ldots,L_{n-1}$ and a global section 
$s_0\in \Gamma(X,L_1)\oplus\cdots\oplus\Gamma(X,L_{n-1})$ whose zero locus is a smooth irreducible curve $C$ which intersects transversally each germ $V_i$ at $p_i$. A generic point $a$ in the associated parameter space 
$$
X^*:=\pp(\Gamma(X,L_1))\times\cdots\times\pp(\Gamma(X,L_{n-1}))
$$
determines a smooth closed curve $C_a$ in $X$, which, for $a$ close enough to the class $a^0\in X^*$ of $s_0$, intersects each germ $V_i$ transversally at a point $p_i(a)$ whose coordinates vary holomorphically with $a$ by the implicit functions theorem. For any function $f$ holomorphic at $p_1,\ldots, p_N$, we define the trace of $f$ on $V:=V_1\cup\cdots\cup V_N$ relatively to $(L_1,\ldots,L_{n-1})$ as the function 
$$
a\longmapsto Tr_V (f)(a):=\sum_{i=1}^N f(p_i(a)),
$$
which is defined and holomorphic for $a$ in an analytic neighborhood of $a^0$.
\vskip 2mm
\noindent

Let us suppose now that $X$ is a toric projective smooth compactification of $U=\cp^n$, endowed with a linear action of an algebraic torus $\mathbb{T}$ that preserves the coordinate hyperplanes $x_i=0,\,\,i=1,\ldots,n$ (see \cite{Dan:gnus}).
Trivially any germ $V_i$ contained in the hypersurface at infinity $X\setminus U$ is algebraic. We can thus suppose that $V$ is contained in $U$ and work with the affine coordinates $x=(x_1,\ldots,x_n)$.
\vskip 2mm
\noindent
Since $U\simeq \cp^n$, its Picard group is trivial and the classes of the irreducible divisors $G_1,\ldots ,G_s$ supported outside  $U$ form a basis for $Pic(X)$. Any globally generated line bundle $L$ on $X$ has thus a unique global section $s_U\in \Gamma(X,L)$ such that $div (s_U)\cap U=\emptyset$. If $s\in \Gamma(X,L)$, the quotient $\frac{s}{s_U}$ defines a rational function without poles on $U\simeq \cp^n$, that is a polynomial in $x$, which gives the local equation for the divisor $H=div(s)$ in the affine chart $U$. Since $L$ is globally generated, a generic section $s\in \Gamma(X,L)$ does not vanish at $0\in U$ and the corresponding polynomial in $x$ has a non-zero constant term. 
\vskip 1mm
\noindent
In our situation of very ample line bundles $L_1,\ldots,L_{n-1}$ on $X$, we can thus use polynomials equations for $C_a$ restricted to the affine chart $U$:
$$
C_a\cap U=\{x=(x_1,\ldots,x_n)\in U,\,\,a_{k0}=P_k(a_k',x),\,\,k=1,\ldots,n-1\},
$$
where $a_k=(a_{k0},a'_k)$ and $P_k(a_k',.)$ are polynomials in $x$ vanishing at $0\in U$. 
\vskip 2mm
\noindent
Since $X$ is toric, we know from \cite{Fult:gnus} that the Chow groups $A_{k}(X)$ are isomorphic to the cohomology groups $H^{2n-2k}(X,\zp)$, for any $k=0,\ldots,n$, and we can identify the Chow group $A_0(X)$ of $0$-cycles with $\zp\simeq H^{2n}(X,\zp)$. We note $[V]$ the class of any closed subvariety $V$ of $X$, $c_1(L)\in H^{2}(X,\zp)$ the first Chern class of any line bundle $L$ on $X$, and we denote by $\smallfrown$ the usual cap product. Our first result is

\begin{thr}
The set $V:=V_1\cup\cdots\cup V_N$ is contained in an algebraic hypersurface $\wt{V}\subset X$ such that
$$
[\wt{V}] \smallfrown \prod_{k=1}^{n-1} c_1(L_k) = N
$$
if and only if for all $i=1,\ldots,n$ the functions $a\mapsto Tr_V (x_i)(a)$ are affine in the constant coefficients $a_0=(a_{10},\ldots,a_{n-1,0})$.
\end{thr}
In the general case, none of the germs $V_i$ has a tangent space at $p_i$ equal to $x_n=0$, in which case the condition ``$Tr_V (x_n)$ affine in $a_0$'' is sufficient for Theorem 2 to hold.
\vskip 1mm
\noindent
If the conditions of Theorem 2 are not satisfied, $V$ can nevertheless be contained in an hypersurface $\wt{V}$ of $X$ such that $[\wt{V}] \smallfrown \prod_{k=1}^{n-1} c_1(L_k) > N$. In that case, traces of affine coordinates are algebraic in $a_0$ and no longer polynomials. Let us mention the following toric Abel-inverse theorem obtained in \cite{WM1:gnus}, chapter 2, as a corollary of Theorem 2, generalizing results of \cite{hp:gnus} and \cite{WM2:gnus}:
\begin{thr}
Let $\phi$ be a holomorphic form of maximal degree on $V$, given by $\phi_i$ on the germ $V_i$, for $i=1,\ldots,N$. There exist an algebraic hypersurface $\wt{V}\subset X$ containing $V$ such that
$[\wt{V}] \smallfrown \prod_{k=1}^{n-1} c_1(L_k) = N$, and a rational form $\Psi$ on $\wt{V}$ such that $\Psi_{|V_i}=\phi_i$ for $i=1,\ldots,N$, if and only if the trace form $Tr_V \phi(a):=\sum_{i=1}^N p_i^*(\phi_i)(a)$ is rational in $a_0$.
\end{thr}
\vskip 2mm
\noindent
Contrary to the projective case handled in \cite{W:gnus}, Theorem 2 does not characterize the class of $\wt{V}$. To do so, we introduce the norm on $V$ relatively to $(L_1,\ldots,L_{n-1})$ of any function $f$ holomorphic at $p_1,\ldots,p_N$,
$$
a\longmapsto N_V (f)(a):=\prod_{i=1}^N f(p_i(a)),
$$
which is defined and holomorphic for $a\in X^*$ close to $a^0$. We then study the degree in $a_0$ of norms of some rational functions on $X$ whose polar divisors generate $Pic_{\qp}(X)$. As in \cite{WM2:gnus}, let us fix very ample effective divisors $E_1,\ldots,E_s$ supported by $X\setminus U$, whose classes form a $\qp$-basis of $Pic_{\qp}(X)$. We can now characterize the class of the interpolating hypersurface~:
\begin{thr}
Suppose that conditions of Theorem $2$ are satisfied. Then the equality $[\wt{V}]=\alpha\in Pic(X)$ holds if and only if
there exist rational functions $f_j\in H^0(X,\mathcal{O}_X(E_j))$ for $j=1,\ldots,s$,  whose norms
$N_V (f_j)$ are polynomials in $a_{10}$ of degree exactly
$$
deg_{a_{10}}\,\,N_V (f_j)= \alpha\cdot [E_j]\smallfrown\prod_{k=2}^{n-1} c_1(L_k)\in \zp_{\ge 0}.
$$
\end{thr}
Note that Bernstein's theorem \cite{Ber:gnus} allows to compute the degrees of intersection in Theorems 2 and 4 as mixed volume of the polytopes associated (up to translation) to the involved line bundles. 
\vskip 1mm
\noindent
If $X=\pp^n$, then $Pic(X)\simeq\zp$ and Theorem $4$ follows from Theorem $2$~:~ if $Tr_V (x_n)$ is affine in $a_0$, then $N_V (x_n)$ has degree $N$ in $a_0$. 
\vskip 2mm
\noindent

The proof of theorem 2 uses a toric generalization of Abel-Jacobi's theorem \cite{Kho:gnus} which gives combinatorial conditions for the vanishing of sums of Grothendieck residues of rational forms in toric varieties, which can be interpreted in term of affine coordinates. 
\vskip 2mm
\noindent
The difficulty to generalize Theorem 2 for other compactifications $X$ of $\cp^n$, as grassmannians or flag varieties, is that there is no natural choice for affine coordinates, so {\it a priori} no grading for the algebra of regular functions over $U=\cp^n$ naturally associated to $X$ (interpolation results in grassmannians would be important for generalizing Theorem 2 to any projective variety $X$ and to union of germs of any dimension $k\le n-1$, by using a grassmannian embedding of $X$ associated to an adequat rank $k$ ample bundle $E$ on $X$). We can hope a generalization to the case of non-projective toric varieties, using blowing-up and essential families of globally generated line bundles, as presented in \cite{WM1:gnus} (chapter 2, section 2). 
\vskip 2mm
\noindent

Section 2 is devoted to the proof of Theorem $2$, and Section 3 to the proof of Theorem $4$. 
\vskip 2mm
\noindent

This article is extracted from my thesis \cite{WM1:gnus}, untitled ``La trace en g\'{e}om\'{e}trie projective et torique'', which is disponible on my home page \\
\texttt{http$:$//www.math.u-bordeaux.fr/$\sim$weimann/}.

\section{Proof of Theorem 2}

\subsection{Direct implication} 

Let us suppose that $V$ is contained in an algebraic hypersurface $\wt{V}$ whose equation in the affine chart $U$ is given by a polynomial $f\in \cp[x_1,\ldots,x_n]$. Since the line bundles $L_1,\ldots,L_{n-1}$ are very ample, the hypothesis on the degree of intersection is equivalent to that for $a$ near $a^0$, the intersection $\wt{V}\cap C_a$ is contained in $U$ and equal to $V\cap C_a$. As explained in the introduction of \cite{WM2:gnus}, the trace of any coordinate $x_i$ on $V$ relatively to $(L_1,\ldots,L_{n-1})$ is equal for $a$ close to $a^0$ to the sum of Grothendieck residues of the rational $n$-form on $U$
$$
\frac{x_i df\land dP_1\cdots\land dP_{n-1}}{f(x)(a_{10}-P_1(a'_1,x))\cdots(a_{n-1,0}-P_{n-1}(a'_{n-1},x))}\,\,,
$$
which we denote, as in \cite{Yger:gnus}, by
$$
Tr_{V}(x_i)(a)={\rm Res}\,
\left[\begin{matrix}
x_i df\land dP_1\cdots\land dP_{n-1}\cr
f,a_{10}-P_1,\ldots  ,a_{n-1,0}-P_{n-1}
\end{matrix}
\right].
$$
Following \cite{cd:gnus} or \cite{WM2:gnus}, the integral representation for the global sum of Grothendieck residues allows us to derivate the trace according to $a_{k0}$ under the integral: 
$$
\partial_{a_{k0}}^{(l)} Tr_{V} (x_i)(a)={\rm Res}\,
\left[\begin{matrix}(-1)^l\,l!\,
x_1\cdots x_i^2 \cdots x_n\frac{df\land dP_1\cdots\land dP_{n-1}}{x_1\cdots x_n}\cr
f,a_{10}-P_1,\ldots,(a_{k0}-P_{k})^{l+1},\ldots,a_{n-1,0}-P_{n-1}
\end{matrix}
\right].
$$
If $h,f_1,\ldots,f_n$ are Laurent polynomials in $t=(t_1,\ldots,t_n)$ with Newton polytopes $P,P_1,\ldots,P_n$, the toric Abel-Jacobi theorem \cite{Kho:gnus} asserts that 
$$
{\rm Res}\,
\left[\begin{matrix}
h \frac{dt_1\cdots\land dt_{n}}{t_1\cdots t_n}\cr
f_1,\ldots,f_n
\end{matrix}
\right]=0
$$ 
as soon as $P$ is contained in the interior of the Minkowski sum $P_1+\cdots +P_n$.
Since $L_k$ is very ample, the support of the polynomial $P_k$ is $n$-dimensional and it is not hard to check that the Newton polytope of the jacobian of the map $(f,P_1,\ldots,P_{n-1})$ translated {\it via} the vector $(1,\ldots,2,\ldots,1)$ (corresponding to  multiplication by $x_1\cdots x_i^2\cdots x_n$) is stricly contained in the Minkowski sum of the Newton polytopes of polynomials $f, a_{10}-P_1,\ldots,a_{n-1,0}-P_{n-1}$ for $l\ge 2$. This shows direct part of Theorem 2.
\begin{rem}
In general, traces of coordinate functions do not depend of $a_0$. If $R_k$ is the unique divisor in $|L_k|$ supported outside $U$, the previous argument yields the implication
$$
h\in H^0(X,\mathcal{O}_X(dR_k))\Rightarrow  deg_{a_{k0}} Tr_V (h) \le d
$$
with equality if the zero set of $h$ has a proper intersection with $X\setminus U$ (which is generically the case since $L_k$ is globally generated). See \cite{WM1:gnus}, Corol. 3.6 p 127.
\end{rem}

\subsection{Converse implication} 

Let us show that $Tr_V (x_i)$ being affine in $a_0$ implies that $Tr_V (x_i^l)$ is polynomial of degree at most $l$ in $a_0$ for any $l\ge 1$. We need an auxiliary lemma generalizing to the toric case the ``Wave-shock equation'' used in \cite{hp:gnus} to show the Abel-inverse theorem. We give a weak version of this lemma, which will be sufficient for our purpose. See \cite{WM1:gnus}, prop. 3.8 p 128, for a stronger version.
\vskip 1mm
\noindent
For $a$ near $a^0$, we use affine coordinates 
$(x_1^{(j)}(a),\ldots,x_{n}^{(j)}(a))$ for the unique point of intersection $p_j(a)$ of $V_j$ with $C_a$. Since $L_k$ is very ample, the monomial $x_i$ occurs in the polynomial $P_k$ with a generically non zero coefficient $a_{ki}$, for $i=1,\ldots,n$.
\begin{lem} 
For  any $i\in\{1,\ldots,n\}$, and any $j\in\{1,\ldots,N\}$, the function $a\mapsto x_i^{(j)}(a)$ (holomorphic at $a^0$) satisfies the following P.D.E:
$$
\partial_{a_{ki}}x_i^{(j)}(a)=x_i^{(j)}\partial_{a_{k0}}x_i^{(j)}(a)
$$
for any $k=1,\ldots,n-1$ and any $a$ close to $a^0$.
\end{lem}
\begin{preuve}
Let us fix $i=1$ for simplicity. Trivially, the equality $a_{k0}=P_k(a_k',x)$ holds for all $k=1,\ldots,n-1$ if and only if $x\in C_a\cap U$, and the complex number
$$
x_1^{(j)}((P_1(a'_1,x),a'_1),\ldots,(P_{n-1}(a'_{n-1},x),a'_{n-1}))
$$
thus represents the $x_1$-coordinate of the unique point of intersection of $V_j$ with the curve $C_a$ passing through $x$.  If $x=(x_1,\ldots,x_n)$ belongs to $V_j$, this complex number, seen as a function of $a'=(a'_1,\ldots,a'_{n-1})$ is thus constant, equal to $x_1$. Differentiating according to the $x_1$-coefficient $a_{k1}$ of $P_k$ gives 
$$
\partial_{a_{k1}}x_1^{(j)}((P_1(a'_1,x),a'_1),\ldots,(P_{n-1}(a'_{n-1},x),a'_{n-1}))
$$
$$
=x_1^{(j)}((P_1(a'_1,x),a'_1),\ldots,(P_{n-1}(a'_{n-1},x),a'_{n-1}))
$$
$$
\times\partial_{a_{k0}}x_1^{(j)}(((P_1(a'_1,x),a'_1),\ldots,(P_{n-1}(a'_{n-1},x),a'_{n-1})).
$$
We can replace $x\in V_j$ with $(x_1^{(j)}(a),\ldots,x_{n}^{(j)}(a))\in V_j$, and the desired relation follows from the equality $P_k(a'_k,(x_1^{(j)}(a),\ldots,x_{n}^{(j)}(a)))=a_{k0}$. $\hfill\square$
\end{preuve}
By induction, this lemma implies the relation 
$$
(l+1)\partial_{a_{ki}}Tr (x_i^l)=l\partial_{a_{k0}}Tr (x_i)^{l+1})
$$
for any $i=1,\ldots,n$, any $k=1,\ldots,n-1$, and all integers $l\in \np$, from wich we easily deduce
$$
{\rm deg}_{a_{k0}} Tr (x_i^l)\le l  \qquad \qquad  (*)
$$
Let $\{f_j=0\}$ be a local (irreducible) equation for the germ $V_j$ and choose a $\cp$-linear combination $u=u_1 x_1+\cdots+u_n x_n
$ of the $x_i$'s such that $\partial_u f(p_j)\ne 0$ for all $j=1,\ldots,N$. We consider then the caracteristic polynomial of $u$:
$$
F_u(Y,a):=\prod_{j=1}^N(Y-u(p_j(a)))
$$
whose coefficients are holomorphic functions near $a^0$. Using Newton formulae relating coefficients of $F_u$ to the trace of the powers of $u$, we deduce from $(*)$ that $F_u$ is polynomial in $a_0=(a_{01},\ldots,a_{0,n-1})$. The function 
$$
Q_{a'}(x):=F_u(u,(P_1(x,a'_1),a'_1),\ldots,(P_{n-1}(x,a'_{n-1}),a'_{n-1}))
$$
is thus a polynomial in $(x_1,\ldots,x_n)$ vanishing on $V_1\cup\ldots\cup V_N$. By hypothesis, we have 
$$
df_j\land dP_1\land\cdots\land dP_{n-1}(p_j)\ne 0\,\,\,{\rm and}\,\,\,\partial_u f(p_j)\ne 0
$$
for all $j=1,\ldots,N$. The implicit function theorem then implies the equivalence
$$
Q_{a'}(x)=0,\,\, x\,\, {\rm near}\,\, C_a\cap U \iff x\in V_1\cup\cdots\cup V_N.
$$ 
Thus, the Zariski closure in $X$ of the algebraic hypersurface $\{Q_{a'}=0\}$ of $U\simeq\cp^n$ does not depend of $a'$ and gives the desired hypersurface $\wt{V}$. $\hfill\square$




\section{Proof of theorem 3}



We can associate to any codimension $2$ closed subvariety $W\subset X$ its dual set $W^*\subset X^*$ associated to the line bundles $(L_1,\ldots,L_{n-1})$, defined by
$$
W^*:=\{a\in X^*,\,\, C_a \cap V\ne \emptyset\}.
$$
From \cite{GKZ:gnus}, this is an hypersurface in the product of projective spaces $X^*$, irreducible if $W$ is,  whose multidegree $(d_1,\ldots,d_{n-1})$ in $X^*$ is given by the intersection numbers
$$
d_i= W\smallfrown \prod_{i=1,i\ne j}^{n-1} c_1(L_i),\,\,\,j=1,\ldots,n-1.
$$
We call the $(L_1,\ldots,L_{n-1})$-resultant of $W$, noted $\mathcal{R}_W$, the multihomogeneous polynomial of multidegree $(d_1,\ldots,d_{n-1})$ vanishing on $W^*$ (it is defined up to a non zero scalar, but this has no consequence for our purpose). By linearity, we generalize this situation to the case of cycles: 
$$
\mathcal{R}_{\sum c_i W_i}:=\prod(\mathcal{R}_{W_i})^{c_i}.
$$
Duality respects rational equivalence: the degree of the resultant of a cycle $W$ only depends of the class of $W$ in the Chow group of $X$ (see \cite{WM1:gnus}, prop. 7 p 100). 
\vskip 2mm
\noindent
From the product formula \cite{PS:gnus}, any rational function $f_j\in H^0(X,\mathcal{O}_X(E_j))$ whose zero divisor $H_j$ intersects properly $\wt{V}$ and $X\setminus U$, gives rise to the equality :
$$
N_{\wt{V}}(f_j)=\frac{\mathcal{R}_{\wt{V}\cdot H_j}}{\mathcal{R}_{\wt{V}\cdot E_j}}.
$$
Since the constant coefficents  $a_0=(a_{10},\ldots,a_{n-1,0})$ do not influence on the asymptotic behavior of the curves $C_a$ outside the affine chart $U$, the resultant $\mathcal{R}_{\wt{V}\cdot E_j}(a)$ does not depend on $a_0$. We thus obtain
$$
{\rm deg}_{a_{10}} N(f_j)={\rm deg}_{a_{10}} \mathcal{R}_{\wt{V}\cdot H_j}\le {\rm deg}_{a_1} \mathcal{R}_{\wt{V}\cdot H_j}={\rm deg}_{a_1}\mathcal{R}_{\wt{V}\cdot E_j}.
$$
Since we deal with homogeneous polynomials in $a_1$, strict inequality in the previous expression is equivalent to the equality
$$
\mathcal{R}_{\wt{V}\cdot E_j}((a_{10},0,\ldots,0),a_2,\ldots,a_{n-1})\equiv 0.
$$
This happens if and only if all subvarieties $C=\{s=0\}$ given by sections $s\in \Gamma(X,\oplus_{k=2}^{n-1}L_{k})$ intersect the set $\wt{V}\cap H_j \cap (X\setminus U)$. By a dimension argument, this would imply that $\wt{V}$ has an irreducible branch contained in $X\setminus U$, which is not possible since $\wt{V}\cap C_a = V\cap C_a \subset U$ for $a$ close to $a^0$.
Thus we have proved the equality:
$$
{\rm deg}_{a_{10}} N_{\wt{V}}(f_j)=[\wt{V}].[E_j]\smallfrown\prod_{k=2}^{n-1} c_1(L_k)
$$
Since the classes $[E_j],\,j=1,\ldots,s$ determine a basis for $A_{n-1}(X)\otimes_{\zp} \qp$, the non degenerated natural pairing between the Chow groups $A_1(X)$ and $A_{n-1}(X)$ shows that hypothesis of theorem $4$ is equivalent to that 
$$
[\wt{V}]\smallfrown \prod_{k=2}^{n-1} c_1(L_k)=\alpha \smallfrown \prod_{k=2}^{n-1} c_1(L_k).
$$
This is equivalent, in turn, to the equality $[\wt{V}]=\alpha$, by Proposition 1.1 in \cite{bg:gnus}, which generalizes the Strong Lefschetz theorem. $\hfill\square$

\end{document}